\documentstyle{czjphys209}
\newcommand{\JJ}{\mathbin{\raisebox{0.25ex}{$
                      \rm\vphantom{I}%
                      \_\hskip -0.25em \_%
                      \vrule width 0.6pt$}}}         %right contraction
\newcommand{\JJB}{\mathop{\JJ \,}\displaylimits_B}
\newcommand{\cl}{C \kern -0.1em \ell}
\newcommand{\openone}{\mbox{1\kern -0.20em I}}
\newcommand{\openC}{\mbox{C\kern -0.55em I\kern 0.15em}}
\newcommand{\openH}{\mbox{I\kern -0.15em H}}
\newcommand{\openL}{\mbox{I\kern -0.15em L}}
\newcommand{\openR}{\mbox{I\kern -0.15em R}}    
\newcommand{\openZ}{\mbox{Z\kern -0.40em Z}}
\newcommand{\Bbbk}{\mbox{I\kern -0.12em k}}
\newcommand{\nn}{\\ \nonumber}
\begin{document}
\title{On the relation of Manin's quantum plane and quantum Clifford
algebras}
\authori{Bertfried Fauser}
\addressi{
Universit\"at Konstanz,
Fachbereich Physik,
Fach M 678\\
78457 Konstanz, Germany\\
E-mail: Bertfried.Fauser@uni-konstanz.de}
\authorii{}     %otherwise {}
\addressii{}
\authoriii{}     %otherwise {}
\addressiii{}
\headtitle{On the relation of Manin's quantum plane\ldots}
\headauthor{Bertfried Fauser}
\specialhead{Bertfried Fauser: On the relation of Manin's quantum plane\ldots}
%%%%%%%%%%%%%% FOR EDITORIAL USE ONLY!!! %%%%%%%%%%%%%%%
\evidence{A}
\daterec{July 24, 2000}    %;-)) final version
\cislo{0}  \year{2000}
\setcounter{page}{1}
\pagesfromto{000--000}
%\makefirsttitle
%%%%%%%%%%%%%%%%%%%%%%%%%%%%%%%%%%%%%%%%%%%%%%%%%%%%%%%%
\maketitle

\begin{abstract}
One particular approach to quantum groups (matrix pseudo groups) provides
the Manin quantum plane. Assuming an appropriate set of non-commuting
variables spanning linearly a representation space one is able to show
that the endomorphisms on that space preserving the non-commutative
structure constitute a quantum group. The non-commutativity of these
variables provide an example of non-commutative geometry.\\   
In some recent work we have shown that quantum Clifford algebras --i.e.
Clifford algebras of an arbitrary bilinear form-- are closely related to
deformed structures as $q$-spin groups, Hecke algebras, $q$-Young
operators and deformed tensor products. The natural question of relating 
Manin's approach to quantum Clifford algebras is addressed here. Some
peculiarities are outlined and explicite computations using the CLIFFORD
Maple package are exhibited. The meaning of non-commutative geometry is
re-examined and interpreted in Clifford algebraic terms.\\[3ex]
15A66 Clifford algebras, spinors\\
11E39 Bilinear and Hermitian forms 
\end{abstract}

\section{Introduction}

The Manin quantum plane is one particular starting point to construct 
quantum groups or matrix pseudo groups \cite{Man,Wor}. The main point in
this construction is, that the non-commutativity of the ``coordinates''
induces a non-commutative tensor product. Since Clifford algebras provide
an ideal tool to describe geometric settings and linear algebra 
\cite{Hes-STS} we try to give an outline in which way the Manin construction 
could be transfered into the Clifford framework. We found in former 
investigations already Hecke algebra representations \cite{Fau-hecke} and
$q$-spin groups \cite{Fau-qgrp} which motivated our work. However, the
common notion of Clifford algebra has to be replaced by that of {\it quantum 
Clifford algebra\/}, that is a Clifford algebra of an arbitrary 
bilinear form, see \cite{FauAbl}. Such type of quantum Clifford algebras 
do no longer possess decomposition theorems, but decompose in general only 
over deformed tensor products. Such deformations need not even to be 
braided. This results also in the fact that the Clifford functor is no 
longer exponential, but might be $q$-exponential in some cases. We examine 
the split case $\cl_{p,p}$ for $p=4$ in the undeformed and deformed case and
show that we can find $x$ and $y$ elements fulfilling the quantum plane 
relations. Also this investigation showed that it is fruitful to
study the decomposition of Clifford algebras to learn about their 
composition.

\section{The Manin quantum plane}

In this section we give a short account on the (two dimensional) quantum
plane as constructed by Manin \cite{Man}. We fix thereby our notation.

Let $V$ be a linear space over a field $\Bbbk$. The dimension of $V$ is given
as the number of {\it linear\/} independent elements. Consider the space
$\otimes^2 V$ of tensors of degree two and impose on elements $x$ and $y$ which
span $V$, i.e. $V = \, <x,y>$, the condition
\begin{eqnarray}
\label{q-plane}
x \otimes y &=& q\, y\otimes x\, .
\end{eqnarray}
Usually the tensor-product is not explicitely written. In a second step, one
considers the endomorphisms of $V$, such that the condition
(\ref{q-plane}) remains invariant ($M \in \mbox{End}(V)$; 
$\vec{x}, \vec{x}^\prime \in V$)
\begin{eqnarray}
M \vec{x} &=& \vec{x}^\prime
\end{eqnarray}
i.e. one requires that
\begin{eqnarray}
x^\prime \otimes y^\prime &=& q\, y^\prime \otimes x^\prime\, .
\end{eqnarray}
The relations imposed on the {\it quantum plane\/} $V$ result in relations of
the matrix elements of $M$ w.r.t. a basis of $V$ and $V^*$. Such matrix
elements are necessarily non-commutative. Obviously one has $V^* =
<\partial_x,\partial_y >$ and the dual quantum plane fulfils the relation
\begin{eqnarray}
\partial_x \otimes \partial_y &=& q^{-1}\, \partial_y \otimes \partial_x.
\end{eqnarray}
This space constitutes the canonical dual space of $V$. Deducing the
corresponding commutation relations one obtains
\begin{eqnarray}
\label{q-com}
M &\cong \,\,
\left( \begin{array}{cc}
a & b \\ c & d 
\end{array} \right) & \nn
ab &=\, q^{-1}\, ba \quad cd    &=\, q^{-1}\, dc \nn
ac &=\, q^{-1}\, ca \quad bd    &=\, q^{-1}\, db \nn
bc &=\, cb          \quad ad-da &=\, (q^{-1}-q)\, bc.
\end{eqnarray}
Finally one can define a trace and a quantum determinant. The later is given as 
$\mbox{Det}_q(M) = ad-q^{-1}bc$ and is a central element of the $GL_q(2)$.

For our purpose it is important to note, that quantum endomorphisms can be
realized as matrices with matrix entries which fulfil some constraints. Such
block matrices constitute then a {\it matrix pseudo group\/} \cite{Wor}, 
since the set of matrices is not a full matrix algebra but 
$\mbox{Mat}_{n\times n}$ mod relations. 

\section{Composition and decomposition}

A main result in the theory of real and complex Clifford algebras is the
occurrence of certain periodicity theorems \cite{Por,Mak,BudTra}. Albert
Crumeyrolle took them as the Mendeljev periodic system of
elementary particles \cite{Mic-Cru}, page X: {\it les th\'eor\`emes 
de p\'eriodicit\'e pour les alg\`ebres de Clifford devraient jouer, pour la 
th\'eorie des particules \'el\'ementaires, un r\^ole analogue \`a celui 
jou\'e par la classification p\'eriodique des \'el\'ements de 
Mendele\"{\i}ev en Chimie\/}. However, such periodicity theorems
break down for quantum Clifford algebras.

\subsection{Periodicity theorems I:}

Let ${\sf H} = (V,Q)$ be a quadratic space, that is a pair of a $\Bbbk$-linear
space $V$ and a quadratic form $Q$ and ${\bf Quad}$ the corresponding 
category. If $\mbox{char~} \Bbbk \not=2$ then $Q$ has a symmetric matrix 
representation. There is a functorial mapping of quadratic spaces called 
$\cl$ into the associative unital algebras
\begin{eqnarray}
{\bf Quad} &\stackrel{\cl}{\longrightarrow} {\bf Alg}.
\end{eqnarray}
The particular map $\gamma : (V,Q) \mapsto \cl(V,Q)$ is called Clifford map
\cite{Che}. The (universal) algebra connected in this way to a quadratic space
is called Clifford algebra. From the basic law, that $V \ni x \mapsto
\gamma_x$ and $\gamma_x\gamma_x = Q(x){\openone}$ one obtains by 
polarization the common anti-commutation relations
\begin{eqnarray}
i) &~~& \gamma_x\gamma_y + \gamma_y\gamma_x = 2g(x,y) \nn
ii) &~~& Q(x+y)-Q(X)-Q(y) =: 2g(x,y).
\end{eqnarray}
During the process of classifying such algebras one obtains periodicity
theorems, as already Clifford proved some of them. In particular one has
\begin{eqnarray}
\label{dec}
i) &~~~ \cl_{p,q+8} &\cong \cl_{p,q} \otimes \openR(16) \nn
ii) &~~ \cl_{p+1,q+1} &\cong \cl_{p,q} \otimes \cl_{1,1} \nn 
                     &&\cong \cl_{p,q} \otimes \openR(2) 
\end{eqnarray}
etc., where $p$ and $q$ are sufficient large to allow the
decomposition. $\cl_{p,q}$ is the real Clifford algebra
$\cl (V,Q)$ with a quadratic form $Q$ of signature $p,q$ and 
$\openR(n)$ is the {\it full\/} $n\times n$ matrix algebra 
over the reals. Obviously one can reduce higher dimensional Clifford algebras
by splitting them into smaller parts. As a consequence of this process one
stays with a few basic irreducible (smallest or elementary) Clifford algebras.
This process can be reversed to {\it compose\/} any real Clifford algebra 
from such building blocks. A particular interesting case for our purpose 
is (\ref{dec}-ii) where we can look at $\cl_{p+1,q+1}$ as a $2 \times 2$ 
matrix algebra over the non-commutative ring $\cl_{p,q}$. Elements 
$M \in \cl_{p+1,q+1}$ can be written as a $2 \times 2$ block matrix 
containing $m_{ij} \in \cl_{p,q}$ as entries, where the indices $i,j$ 
are w.r.t. $\openR(2)$. At this stage it is important that the 
classification of Clifford algebras is identical to the classification 
of the involved quadratic forms $Q$.

We can learn from the decomposition theorems moreover, that the Clifford
functor $\cl$ is exponential. That is a direct sum of two spaces is mapped onto
the tensor-product of Clifford algebras
\begin{eqnarray}
\cl(U\oplus V,Q_1 \perp Q_2) &\cong \cl(U,Q_1) \otimes \cl(V,Q_2).
\end{eqnarray}

\subsection{Quantum Clifford algebras}

Let ${\sf R} =\, (V,B)$ be a reflexive space, that is a $\Bbbk$-linear 
space $V$ and a bilinear not necessarily symmetric form 
$B : V \times V \mapsto \Bbbk$. To define the associated 
{\it quantum Clifford algebra\/} we start defining the exterior algebra 
$\bigwedge V$ over $V$,
\begin{eqnarray}
\bigwedge V &=&
\Bbbk \oplus V \oplus \wedge^2 V \oplus \ldots 
\oplus \wedge^n V \oplus \ldots \,.
\end{eqnarray}
We introduce the dual space $V^*$ of linear forms on $V$ and a so-called
Euclidean dual isomorphism $\ast : V \mapsto V^*$ to introduce a dual 
basis as a particular isomorphism $V^* \ni i_x \cong x^*$. Let $x,y \in V$ 
and $x^* \in V^*$ the image of $x$ in $V^*$. The contraction is defined as
\begin{eqnarray}
i_x(y) \,=\, x^*(y) &=& x \JJB y \,=\, B(x,y).
\end{eqnarray}
We write $\JJ$ for short if it is clear which bilinear form is involved 
in the definition of the contraction. Following the ideas of Chevalley 
\cite{Che} we note that with $x \in V$
\begin{eqnarray}
\gamma_x &=& x\JJB\,\,+x \wedge
\end{eqnarray}
is a Clifford map of ${\sf R} = \,(V,B)$ into ${\bf Alg}$. The particular
Clifford algebra appears as strict subalgebra of the endomorphism algebra 
$\cl(V,B) \,\subset\, \mbox{End}(\bigwedge V)$ of 
the exterior algebra $\bigwedge V$. The action of $\gamma_x$ --we use 
$x$ for short from now on-- induces an action on the whole exterior 
algebra demanding that ($x,y \in V$; $u,v,w \in \bigwedge V$)
\begin{eqnarray}
\label{con}
i)	&~~ x \JJB y             &=\, B(x,y) \nn
ii)	&~~ x \JJB (u \wedge v)  &=\, (x \JJB u) \wedge v
 	    + \hat{v} \wedge (x \JJB u) \nn
iii)	&~~ (u \wedge v) \JJB w  &=\, u \JJB (v \JJB w),
\end{eqnarray} 
where $\hat{v}$ is the Grassmann grade involution defined as 
$\hat{v} = (-1)^{\partial v} v$, and $\partial v$ is the Grassmann grade
of the element $v$. See \cite{Fau-wick} for an Hopf algebraic 
explanation of the grade involution as the antipode.

Quantum Clifford algebras are denoted as $\cl(V,B)$ they have been 
investigated in \cite{Fau-hecke,Fau-trans,FauAbl}. As a matter of fact 
we lack a classification of this type of algebras.

\subsection{Periodicity theorems II:}

We have the following theorem \cite{FauAbl}:\\
Let $\cl(V,Q)$ and $\cl(V,B)$ be two Clifford algebras over the space 
$V$ such that $Q(x) = B(x,x)$, i.e. $Q$ and $B$ describe the same quadratic
form, then there exists a unique $\openZ_2$-graded Wick-isomorphism $\phi$
connecting the algebras via
\begin{eqnarray}
\phi \circ \cl(V,Q) &=& \cl(V,B).
\end{eqnarray}

This theorem can be used to induce periodicity theorems on $\cl(V,B)$. 
However, this can be achieved {\it only\/} if we allow {\it deformed\/}
tensor products in general. Let $W=V \perp_Q U = V \oplus U$ furthermore 
$Q_1:=Q\vert_V$, $Q_2 = Q\vert_U$ the restrictions of $Q$ on $V$ and $U$, 
then we decompose $\cl(V,B)$ as
\begin{eqnarray}
\cl(W,B) &=& \phi \circ \cl(W,Q) \,=\,
\phi \circ \cl(V\oplus U, Q_1 \oplus Q_2) \nn
&=& \phi \circ (\cl(V,Q_1) \otimes \cl(V,Q_2)) \nn
&=&(\phi\vert_V \circ \cl(V,Q_1))
   \otimes\!\!\vert_\phi
   (\phi\vert_U \circ \cl(U,Q_2)) \nn
&=&\cl(V,B_1) \otimes\!\!\vert_\phi \cl(U,B_2). 
\end{eqnarray}
But in general $\phi$ does {\it not\/} factorize w.r.t. $Q$, as indicated by 
the $\phi$-deformed tensor product. We note that the Clifford functor
is {\it not\/} exponential for quantum Clifford algebras in general.
There might be a chance to find $\phi$-deformed exponentials. However,
the above defined deformation is not e.g. braided by construction but the
most general possible deformation. This follows from the observation that
we have no further freedom in the reflexive space ${\sf R}$ and that
the Clifford functor is injective.

\section{Reconstruction of the quantum plane}

\subsection{The matrix window}

We know from the structure theory of real Clifford algebras that such 
algebras possess faithful irreducible representations on spinor spaces 
over distinguished rings. Let 
$\openL \in \{\openR, \openC, \openH, {}^2\openR, {}^2\openH \}$,
then one obtains representations $\rho$ with
\begin{eqnarray}
\rho(\cl_{p,q}) &\cong& \mbox{Mat}_{\openL}(n\times n),
\end{eqnarray}
where $\openL$ is chosen accordingly to the Radon-Hurwitz number. 
It is important to note that such representations map Clifford algebras
onto {\it full\/} matrix algebras and that there are no further relations 
in force.

Using tensor products of matrices, one has to give an identification 
of bases of tensor factors. This means that one has no {\it a priori\/} 
information about the relation of bases of different tensor factors. 
This is especially important if one wants to switch entries in a tensor 
product, e.g. ($x_{\{e\} }$ is $x$ in basis $\{ e\}$)
\begin{eqnarray}
sw(x_{\{e\}} \otimes y_{\{f\}}) 
  &=& 
y_{\{e\}} \otimes x_{\{f\}}.
\end{eqnarray}
The Clifford approach is intrinsically invariant and basis free. 
Furthermore, if one wants to {\it compose\/} Clifford algebras, one 
can learn about this construction by {\it decomposing\/} larger 
Clifford algebras into factors.

It is important to note that only spinor bases provide representations 
of Clifford algebras as full matrix algebras. Representation indices 
are spinor indices.

Spinor spaces can be constructed inside the Clifford algebra as minimal 
left (right) ideals generated from primitive (indecomposable) idempotent 
elements $f_{ii}$. One finds
\begin{eqnarray}
{\cal S} &\cong& \cl f_{11} \nn
\{ e \} & =& \{f_{11},x_2 f_{11},\ldots,x_n f_{11} \}
\,=\, \{ f_{11},f_{21},\ldots,f_{n1} \}
\end{eqnarray}
where the spinor basis $\{ e\}$ is generated by $f_{11}$ and the $x_i$ 
which intertwine the different primitive idempotent elements of the 
Clifford algebra
\begin{eqnarray}
x_r f_{11} &=& f_{rr} x_r
\end{eqnarray}
in a compatible way, i.e. $x_r = x_s x_t$ if $s,t <r$ if possible. Such 
bases can be found e.g. in \cite{Fau-aachen}.

\subsection{Computational examples}

\subsubsection{Undeformed case:}

The space time split \cite{Wey} as the modulo $(1,1)$ periodicity 
\cite{Mak} --called also conformal split w.r.t. $\cl_{1,1}$-- have 
been extensively studied in \cite{Hes-STS}. We restrict 
ourself to the split Clifford algebras $\cl_{p,p}$ and $p=2$. From 
(\ref{dec}) we have
\begin{eqnarray}
\cl_{2,2} &\cong& \cl_{1,1} \otimes \cl_{1,1} 
\,\cong\, \cl_{1,1} \otimes \openR(2) \nn
\cl_{2,2} &\ni& 
M = \left(\begin{array}{cc} a & b \\ c & d \end{array}\right)
\quad\mbox{with~~} a,b,c,d \in \cl_{1,1}\,.
\end{eqnarray}
Let $\{e_i\}$ and $\{f_j\}$ be the generators of the first and second
$\cl_{1,1}$ algebras. Both of them are isomorphic to $\openR(2)$ and 
possess matrix bases $E_{ij}$ and $F_{ij}$. In a spinor basis $\xi_\alpha$
this is $(E_{ij})_{\alpha\alpha^\prime} = 
\delta_{i\alpha}\delta_{j\alpha^\prime}$. Any element of $\cl_{2,2}$ 
can be written as
\begin{eqnarray}
M &=& \left(\begin{array}{cc} a & b \\ c & d \end{array}\right)
\,=\, a F_{11}+ b F_{12}+ c F_{21}+ d F_{22}. 
\end{eqnarray}
It can be checked, that $\{e_i \otimes f_1 f_2\}$ and 
$\{\openone \otimes f_j\}$ generate $\cl_{2,2}$. 

The {\bf spin} and {\bf pin} groups of Clifford algebras are defined to
be sub groups of the Clifford Lipschitz group $\Gamma(p,q)$, which is 
the subgroup of multiplicative units of the Clifford algebra, denoted
as $\cl^*$,  preserving the $\openZ_n$-multi-vector grading 
\begin{eqnarray}
\Gamma(p,q) &:=& \{ g \in \cl_{p,q}^* \,:\, g V g^{-1} \in V \}. 
\end{eqnarray} 
Using the reversion ($(AB)\tilde{~}=A\tilde{~}B\tilde{~}$,
$e_i\tilde{~}=e_i$, $\openone\tilde{~}=\openone$) and the grade involution
$\hat{~} : V \mapsto -V$ as defined in (\ref{con}) one composes two 
anti-involutions $\alpha_\epsilon$ as $\alpha_1 = \tilde{~}$ and 
$\alpha_{-1} = \hat{~} \circ \tilde{~}$. Define 
\begin{eqnarray}
\mbox{\bf pin}(p,q) &:=& \{g \in \Gamma(p,q)\,:\, g\alpha_{-1}(g)=\pm 1 \} \nn
\mbox{\bf spin}(p,q) &:=& \mbox{\bf pin} \cap \cl^+(p,q),  
\end{eqnarray} 
where $\cl^+$ is the even part of $\cl$ w.r.t. the $\openZ_2$-grading. 
The elements of the Clifford--Lipschitz group are called {\it versors\/} 
and are decomposable into Clifford products of vectors --i.e. grade 1-- 
elements of $\cl$.

Conformal transformations have been investigated by Vahlen \cite{Vah}.
In Clifford terms one is interested in elements $M$ from 
$\cl_{1,1}\otimes \openR(2)$ which respect the $\openZ_4$-grade of 
$\cl_{2,2}$. This imposes certain restrictions on the matrix elements of
$M$. One obtains \cite{Mak}:
\begin{eqnarray}
i) &~~& aa\tilde{~}, bb\tilde{~}, cc\tilde{~}, dd\tilde{~} \,\in\,\openR \nn
ii) &~~& ac\tilde{~}, bd\tilde{~} \,\in\, V \nn
iii) &~~& avb\tilde{~} + bv\tilde{~} a\tilde{~},
          cvd\tilde{~} + dv\tilde{~} c\tilde{~} \,\in\, \openR 
\quad\forall \,v\,\in\, V \nn
iv) &~~& avd\tilde{~} + bv\tilde{~} c\tilde{~} \,\in\, V 
\quad\forall \,v\,\in\, V \nn
v) &~~& a\hat{b}\tilde{~} = b\hat{a}\tilde{~},\quad 
        c\hat{d}\tilde{~} = d\hat{c}\tilde{~} \nn
vi) &~~&  a\hat{d}\tilde{~} - b\hat{c}\tilde{~} \,=\, \pm 1. 
\end{eqnarray}
The last two conditions are normalization conditions and $vi)$ could be
called pseudo determinant. Compare this with the relations (\ref{q-com}) 
and the $q$-determinant there. This type of matrices is suitable to model
conformal geometry.

\subsubsection{Deformed case:}

We use CLIFFORD, a Maple V package developed by Rafa{\l} Ab{\l}amowicz
\cite{Abl} to demonstrate that we can indeed find nontrivial elements 
fulfilling the Manin quantum plane relations. Therefore we search for 
two elements $x$ and $y$ in $\cl(V,B)$, where $B$ is defined to be
\begin{eqnarray}
B &:=& 
\left( \begin{array}{cc}
g_{11} & g_{12}+A_{12} \\ g_{12}-A_{12} & g_{22} \end{array} \right).
\end{eqnarray}
A spinor basis consists of an idempotent element and a second linear
independent element. Starting with arbitrary elements, CLIFFORD finds
after having set an arbitrary parameter to zero (giving results only, 
because lack of space)
\begin{eqnarray}
x &:=& \frac{1}{2}(1+A_{12})\openone 
     + \frac{\sqrt{g_{11}(\mbox{det}(g)+1)}}{2g_{11}}{\bf e}_1
     + \frac{1}{2}{\bf e}_1 \wedge {\bf e}_2
\nn
y &:=&  \frac{1}{2}\sqrt{g_{11}(\mbox{det}(g)+1)}\openone
       +\frac{1}{2}(1-g_{12}){\bf e}_2
       +\frac{1}{2}g_{11} {\bf e}_1, 
\end{eqnarray}
where $\mbox{det}(g)$ is the determinant of the symmetric part of $B$.
After solving the $RXY = YXR$ equation, we find two linearly independent
solutions for $R$ 
\begin{eqnarray}
R_1 &:=&  (A_{12}-1)\openone
         +\frac{\sqrt{(\mbox{det}(g)+1)g_{11}}}{g_{11}} {\bf e}_1
         +{\bf e}_1 \wedge {\bf e}_2
\nn
R_2 &:=&  \frac{\sqrt{(\mbox{det}(g)+1)g_{11}}}{g_{11}}\openone
         +\frac{1-g_{12}}{g_{11}} {\bf e}_1
         +{\bf e}_2. 
\end{eqnarray}
Obviously we could also find $q$, however as an algebra element,
which relates the two basis elements $x$ and $y$. From this formulas, 
it is clear that the antisymmetric part is essential to be able to
calculate such $R_i$ elements. The term quantum Clifford algebra was
chosen to exhibit this connection to deformed structures in the case
of non-symmetric bilinear forms.

The above discussed Vahlen matrices of conformal transformations
will be affected by such a deformation too. One would expect to
find as a special case $q$-deformed conformal transformations. 

In Hahn \cite{Hah} it was shown that Clifford algebras have a 
separability idempotent element which allows one to inject the 
algebra, here $\cl_{1,1}$ into its enveloping algebra, here 
$\cl_{2,2} \cong \cl_{1,1} \otimes \cl_{1,1}$. The former algebra 
is a bi-module under the action of the enveloping algebra and the 
unit is the inverse image of the separability idempotent. This 
structure defines in fact the (undeformed) co-product. The deformed
case of this structure will be investigated in detail elsewhere.  

\section{conclusion}

We showed that quantum Clifford algebras can be used to investigate
deformed tensor products and the groups acting on modules over such
deformed spaces. Examining the decomposition of ordinary and quantum 
Clifford algebras, we have been able to show that quantum plane 
relations can be modelled in this framework. This supports our finding
of $q$-spin groups and Hecke algebra representations in previous 
works. Vahlen matrices and the conformal split seem to be the 
undeformed counterparts of the quantum plane. This sheds some light 
on the geometrical meaning of quantum planes and non-commutative geometry,
which in our opinion is a symmetry of ``composed'' objects. In geometrical 
terms we have introduced inhomogeneous coordinates and a metric by 
introducing a Cayley-Klein type measure induced by the factorization 
\cite{Wey}. An important issue for further research is the connection of
quantum Clifford algebras and Hopf algebras as examined in 
\cite{Fau-wick}. A more satisfying theory would definitively incorporate
the Hopf algebra structure and develop all groups and invariants from
that starting point.

\end{document}